\documentclass[11pt]{article}
\usepackage{amsfonts}
\newtheorem{Theoreme}{Th\'eor\`eme}[section]
\newtheorem{Theorem}{Theorem}[section]

\newtheorem{Defi}[Theorem]{Definition}

\newtheorem{Lemma}[Theorem]{Lemma}

\newtheorem{Remarque1}[Theoreme]{Remarque}

\newtheorem{Remark1}[Theorem]{Remark}

\newtheorem{Propo}[Theorem]{Proposition}

\newenvironment{Proof}{\medbreak{\noindent\bf Proof }}{~{\hskip3pt$\bullet$\bigbreak}}
\newenvironment{Aknow}{\medbreak{\noindent\bf Aknowledgement: }}{~{\bigbreak}}

\newenvironment{Remark}{\begin{Remark1}\em}{\end{Remark1}}

\renewcommand{\Im}{{\cal F}}

\newcommand{\hp}{\hskip 3pt}
\newcommand{\hph}{\hskip 8pt}

\newcommand{\K}{{\mathbb K}}

\newcommand{\N}{{\mathbb N}}
\newcommand{\Esp}{{\mathbb E}}

\newcommand{\Ha}{{\mathcal H}}
\newcommand{\Pa}{{\mathcal T}}

\newcommand{\ds}{\displaystyle}
\newcommand{\un}{{\mathbf 1}}

\newcommand{\Dim}{\mbox{Dim}}

\begin{document}
\title{On entropy and Hausdorff dimension of measures defined through a non-homogeneous Markov process} 
\author{Athanasios BATAKIS}
\date{}
\maketitle
\begin{abstract}
In this work we study the Hausdorff dimension of measures whose weight distribution satisfies a markov non-homogeneous property. We prove, in particular,  that the Hausdorff dimensions of this kind of measures coincide with their lower R\'enyi dimensions (entropy). Moreover, we show that the Tricot dimensions (packing dimension) equal the upper R\'enyi dimensions. 

As an application we get a continuity property of the Hausdorff dimension of the measures, when it is seen as a function of the distributed weights under the $\ell^{\infty}$ norm.
\end{abstract}

\section{Introduction}

Let us consider the dyadic tree (even though all the results in this paper can be easily generalised to any $\ell$-adic structure, 
$\ell\in\N$), let $\K$ be its limit (Cantor) set and note $\left({\mathcal F}_n\right)_{n\in\N}$ the  associated filtration with the usual $0-1$ encoding. We are interested in Borel measures $\mu$ on $\K$ constructed in the following way: 
Take $(p_n,q_n)_{n\in\N}$ a sequence of couples of real numbers satisfying $0\le p_n,q_n\le 1$.  

Let $I=I_{\epsilon_1,...,\epsilon_n}$ be a cylinder of the $n$th generation and $IJ=I_{\epsilon_1,...,\epsilon_n,\epsilon_{n+1}}$ a subcylinder of the $(n+1)$th generation, where $\epsilon_1,...,\epsilon_n,\epsilon_{n+1}\in\{0,1\}$. 
\nocite{Pey} \nocite{Wu98}
The mass distribution of $\mu_{|I}$ will be as follows: $\mu(I_0)=p_0\; ,\; \mu(I_1)=1-p_0$ et  
\begin{equation}\label{markov}
\frac{\mu(IJ)}{\mu(I)}=\Big\{^{\displaystyle p_n\un_{\{\epsilon_{n+1}=0\}}+(1-p_n)\un_{\{\epsilon_{n+1}=1\}}\mbox{ , if }\epsilon_n=0}_{\displaystyle q_n\un_{\{\epsilon_{n+1}=0\}}+(1-q_n)\un_{\{\epsilon_{n+1}=1\}}\mbox{ , if }\epsilon_n=1}
\end{equation}
We use the notation $\dim_{\Ha}$ for the Hausdorff dimension and $\dim_{\Pa}$ for the packing (Tricot) dimension.
\begin{Defi}
If $\mu$ is a measure on $\K$, we will denote by $h_*(\mu)$ the lower entropy of the  measure :
$$h_*(\mu)=\liminf_{n\to\infty}\frac{-1}{n} \sum_{I\in{\mathcal F}_n}\log\mu(I)\cdot\mu(I),$$ 
by 
$h^*(\mu)$ the upper entropy of the  measure :
$$h^*(\mu)=\limsup_{n\to\infty}\frac{-1}{n} \sum_{I\in{\mathcal F}_n}\log\mu(I)\cdot\mu(I),$$ 
by $\dim_*(\mu)$ the lower Hausdorff dimension of $\mu$:
$$\dim_*{\mu}=\inf\{\dim_{\Ha} E\; ;\; E\subset\K\mbox{ and }\mu(E)>0\}$$ and by 
$\dim^*(\mu)$ the upper Hausdorff dimension of $\mu$:
$$\dim^*{\mu}=\inf\{\dim_{\Ha} E\; ;\; E\subset\K\mbox{ and }\mu(\K\setminus E)=0\}.$$
In the same way we define
the lower packing dimension (Tricot dimension) of $\mu$:
{\em $$\Dim_*{\mu}=\inf\{\dim_{\Pa} E\; ;\; E\subset\K\mbox{ and }\mu(E)>0\}$$} and by 
{\em $\Dim^*(\mu)$} the upper Hausdorff dimension of $\mu$:
{\em $$\Dim^*{\mu}=\inf\{\dim_{\Pa} E\; ;\; E\subset\K\mbox{ and }\mu(\K\setminus E)=0\}.$$}
\end{Defi}
One can show that (see \cite{Batakis5},\cite{BH})
$$\dim_*(\mu)\le h_*(\mu)\le h^*(\mu)\le \Dim^*(\mu),$$
and there are examples of these inequalities being strict, even when the measure $\mu$ is rather ``regular''.

It is also well known (cf \cite{Fal97}, \cite{Bil}, \cite{Mattila}, \cite{Fan}, \cite{You}, \cite{Renyi} and \cite{Heurteaux}) that 
$$\displaystyle \dim_*(\mu)=\inf\mbox{ess}_{\mu}\liminf_{n\to\infty}\frac{\log\mu(I_n(x))}{-n\log 2}$$ 
and 
$$\displaystyle \dim^*(\mu)=\sup\mbox{ess}_{\mu}\liminf_{n\to\infty}\frac{\log\mu(I_n(x))}{-n\log 2},$$ 
where $I_n(x)$ is the dyadic cylinder of the $n$th generation containing $x$, $\inf\mbox{ess}_{\mu}$ is the essential infimum and $\sup\mbox{ess}_{\mu}$ is the essential supremum,  taken over $\mu$-almost all $x\in\K$.

In the case of measures defined by (\ref{markov}) we can use tools developed in \cite{Batakis} and \cite{Bata2} to prove they are exact, i.e. that $\dim_*(\mu)= \dim^*(\mu)$ or equivalently that
$$\displaystyle\liminf_{n\to\infty}\frac{\log\mu(I_n(x))}{-n\log 2}=\dim_*(\mu),\mbox{ for }\mu\mbox{-almost all }x\in\K$$ and therefore $\dim_*(\mu)=\dim^*(\mu)$. 
However, theorem \ref{letsgo1} implies this statement. 

In general, there is no trivial inequality relation between $h_*(\mu)$ and $\dim^*(\mu)$. Furthermore, it is easy to construct measures $\mu$ satisfying (\ref{markov}) such that $h_*(\mu)\not= h^*(\mu)$ which shows that the sequence of functions $\ds \frac{\log\mu(I_n(x))}{-n\log 2}$ does not necessarily converge (in any space). 

The proof of theorem \ref{letsgo1} implies that there is a sequence $(c_n)_{n\in\N}$ of real numbers such that 
$$\lim_{n\to\infty}\left[\frac{\log\mu(I_n(x))}{-n\log 2}-c_n\right]=0,$$ where $\ds c_n=\frac{-1}{n\log 2}\sum_{I\in{\mathcal F}_n}\log(\mu(I))\mu(I)$. This can be seen as a Shannon-McMillan-type theorem generalised to measures defined through non-homogeneous Markov chains. 

Remark that the tools of \cite{KaP} and \cite{Ka} can be applied to give the same results for {\em ``almost every''} measure $\mu$ satisfying (\ref{markov}). Other results in this sense involving coloring of graphs are proposed in \cite{Fathi}. 

A. Bisbas and C. Karanikas \cite{Biskar2} have already partially proved the conclusions of theorem \ref{letsgo1}, for this kind of measures,  under some assumptions on the sequences $(p_n,q_n)_{n\in\N}$. In particular they prove the theorem when the  sequences $(p_n,q_n)_{n\in\N}$ are uniformly bounded away from $0$ and $1$, which is the case of a perturbation of an homogeneous Markov chain. We thank A. Bisbas for communicating to us this article.

\begin{Theorem}\label{letsgo1}
If $\mu$ satisfies (\ref{markov}) then 
{\em$$\dim_*(\mu)=\dim^*(\mu)=h_*(\mu)\mbox{ and }\Dim_*(\mu)=\Dim^*(\mu)=h^*(\mu).$$}
\end{Theorem}
Using the same type of arguments we also obtain the following continuity result.
\begin{Theorem}\label{letsgo2}
Let $\mu$ and $\mu'$ be measures defined by (\ref{markov}) and the corresponding sequences $(p_n,q_n)_{n\in\N}$ and $(p_n',q_n')_{n\in\N}$ respectively. Then $ |\dim_*(\mu)-\dim_*(\mu')|$ and {\em $|\Dim_*(\mu)-\Dim_*(\mu')| $} go to $0$ as $||(p_n,q_n)_{n\in\N}-(p_n',q_n')_{n\in\N}||_{\infty}$ tends to $0$. 
\end{Theorem}

\section{Lemmas and preliminary results}
Let us introduce some notation: 
for $p\in[0,1]$ we note 
$$h(p)=p\log p+(1-p)\log(1-p)$$ and if $I=I_{\epsilon_1,...,\epsilon_{n-1}}\in{\mathcal F}_n,$ 
let us also set $$\gamma(I,n)=\sum_{i=0,1}\log\left(\frac{\mu(II_i)}{\mu(I)}\right)\frac{\mu(II_i)}{\mu(I)}.$$ 
Remark that for $n\in\N$ and $I\in{\mathcal F}_{n-1}$,$\gamma(I,n)$ is equal to
$h(p_n)$ if $\epsilon_{n-1}=0$ and to $h(q_n)$ if $\epsilon_{n-1}=1$  and therefore $\gamma(I,n)$ is absolutely bounded by $\log 2$.

Let us start with the following easy lemma. 
\begin{Lemma}\label{first}
For all $n,k\in\N$ and all $I\in{\mathcal F}_{n-1}$ we can write  
\begin{equation}\label{easy}
\sum_{K\in{\mathcal F}_k}\log\left(\frac{\mu(IK)}{\mu(I)}\right)\frac{\mu(IK)}{\mu(I)}= \gamma(I,n)+
\sum_{i=0,1}\frac{\mu(II_i)}{\mu(I)}\sum_{K\in{\mathcal F}_{k-1}}\log\left(\frac{\mu(II_iK)}{\mu(II_i)}\right)\frac{\mu(II_iK)}{\mu(II_i)}.
\end{equation}
where $I_0$ and $I_1$ are the two cylinders of the first generation.

Furthermore, if we note $a_n^k(I)$ and $b_n^k(I)$ respectively the quantities 
$$a_n^k(I)=\sum_{K\in{\mathcal F}_{k-1}}\log\left(\frac{\mu(II_0K)}{\mu(II_0)}\right)\frac{\mu(II_0K)}{\mu(II_0)} \mbox{ and }$$
$$b_n^k(I)=\sum_{K\in{\mathcal F}_{k-1}}\log\left(\frac{\mu(II_1K)}{\mu(II_1)}\right)\frac{\mu(II_1K)}{\mu(II_1)}$$
then $a_n^k(I)=a_n^k(I')$ and $b_n^k(I)=b_n^k(I')$, for all $I,I'\in{\mathcal F}_n$.
\end{Lemma}
\begin{Proof}
We have 
\begin{eqnarray}\label{easy1}
&&\sum_{K\in{\mathcal F}_k}\log\left(\frac{\mu(IK)}{\mu(I)}\right)\frac{\mu(IK)}{\mu(I)}= \cr
&&\sum_{i=0,1}\sum_{K\in{\mathcal F}_{k-1}}\log\left(\frac{\mu(II_iK)}{\mu(I)}\right)\frac{\mu(II_iK)}{\mu(I)}= \cr
&&\sum_{i=0,1}\sum_{K\in{\mathcal F}_{k-1}} \log\left(\frac{\mu(II_iK)}{\mu(II_i)}\right)\frac{\mu(II_iK)}{\mu(I)}+ \sum_{i=0,1}\log\left(\frac{\mu(II_i)}{\mu(I)}\right)\frac{\mu(II_i)}{\mu(I)}
\end{eqnarray}
Since we have set
$$\gamma(I,n)= \sum_{i=0,1}\log\left(\frac{\mu(II_i)}{\mu(I)}\right)\frac{\mu(II_i)}{\mu(I)},$$
the equalities (\ref{easy1}) give
$$\sum_{K\in{\mathcal F}_k}\log\left(\frac{\mu(IK)}{\mu(I)}\right)\frac{\mu(IK)}{\mu(I)}= 
\gamma(I,n)+ \sum_{i=0,1}\frac{\mu(II_i)}{\mu(I)}\sum_{K\in{\mathcal F}_{k-1}}\log\left(\frac{\mu(II_iK)}{\mu(II_i)}\right)\frac{\mu(II_iK)}{\mu(II_i)}.$$
Is is immediate that $0\le-\gamma(I,n)\le\log 2$.
By the construction of the measure, the quantities $a_n^k(I)$ and $b_n^k(I)$ do not depend on the cylinder $I$ but only on the cylinder's generation $n$ and this ends the proof.
\end{Proof}

\begin{Remark}
Since the quantities $a_n^k(I)$ and $b_n^k(I)$ depend only on the generation of $I$ and on $k$, we can note $a_n^k=a_n^k(I)$ and $b_n^k=b_n^k(I)$  for $I\in{\mathcal F}_n$. We also note $\ds\Delta_n^k=\frac{1}{k}|a_n^k-b_n^k|$.
\end{Remark}

We also need the following technical estimates.
\begin{Lemma}\label{MA401}
For all $p,q\in[0,1]$ we have $\ds |h(p)-h(q)|\le(1-|p-q|)\log 2$. Furthermore, for all $k\in\N$ and all $\alpha>0$,
$$\frac{(1- |p-q|)\log 2}{k+1}+|p-q| \left(1-\frac{1}{k+1}\right)\alpha\le \max\left\{\left(1-\frac{1}{k+1}\right)\alpha,\frac{\log2}{k+1}\right\}.$$
\end{Lemma}
The proof uses elementary 2-dimensional calculus and is therefore omitted. 
\begin{Propo}\label{mainpro}
Let $I,I'$ be two cylinders of the $n$th generation. Then
$$\frac{1}{k}\left| \sum_{K\in{\mathcal F}_k}\log\left(\frac{\mu(IK)}{\mu(I)}\right)\frac{\mu(IK)}{\mu(I)} - \sum_{K\in{\mathcal F}_k}\log\left(\frac{\mu(I'K)}{\mu(I')}\right)\frac{\mu(I'K)}{\mu(I')}\right|<\eta(k)$$
where $\eta$ is a positive function, not depending on $n$, such that $\eta(k)$ goes to $0$ as $k$ tends to $\infty$.
\end{Propo}

\begin{Proof}
Take any two cylinders $I=I_{\epsilon_1,...\epsilon_n},I'=I_{\epsilon_1',...\epsilon_n'}$ of the $n$th generation. If $\epsilon_n=\epsilon_n'$ then by definition of the measure $\mu$ we get 
$$\frac{1}{k}\left| \sum_{K\in{\mathcal F}_k}\log\left(\frac{\mu(IK)}{\mu(I)}\right)\frac{\mu(IK)}{\mu(I)} - \sum_{K\in{\mathcal F}_k}\log\left(\frac{\mu(I'K)}{\mu(I')}\right)\frac{\mu(I'K)}{\mu(I')}\right|=0.$$
If $\epsilon_n\not=\epsilon_n'$, using lemma \ref{MA401}, lemma \ref{first} and the notation therein we obtain:
\begin{eqnarray}\label{stage1}
 \Delta_{n-1}^{k+1}&=&\left| \frac{1}{k+1} \sum_{K\in{\mathcal F}_{k+1}}\log\left(\frac{\mu(IK)}{\mu(I)}\right)\frac{\mu(IK)}{\mu(I)} - \frac{1}{k+1}\sum_{K\in{\mathcal F}_{k+1}}\log\left(\frac{\mu(I'K)}{\mu(I')}\right)\frac{\mu(I'K)}{\mu(I')}\right|=\cr
&=& \Bigg| \frac{\gamma(I,n)-\gamma(I',n)}{k+1}+ \frac{1}{k+1}  \frac{\mu(II_0)}{\mu(I)}\sum_{K\in{\mathcal F}_k}\log\left(\frac{\mu(II_0K)}{\mu(II_0)}\right)\frac{\mu(II_0K)}{\mu(II_0)}\; +\cr
&+&\frac{1}{k+1}  \frac{\mu(II_1)}{\mu(I)}\sum_{K\in{\mathcal F}_k}\log\left(\frac{\mu(II_1K)}{\mu(II_1)}\right)\frac{\mu(II_1K)}{\mu(II_1)} 
\;-\cr
&-& \frac{1}{k+1}  \frac{\mu(I'I_0)}{\mu(I')}\sum_{K\in{\mathcal F}_k}\log\left(\frac{\mu(I'I_0K)}{\mu(I'I_0)}\right)\frac{\mu(I'I_0K)}{\mu(I'I_0)}\; -\cr
&-& \frac{1}{k+1}  \frac{\mu(I'I_1)}{\mu(I')}\sum_{K\in{\mathcal F}_k}\log\left(\frac{\mu(I'I_1K)}{\mu(I'I_1)}\right)\frac{\mu(I'I_1K)}{\mu(I'I_1)}
\Bigg|=\cr
&=&\left| \frac{ h(p_n)-h(q_n) }{k+1}+ \frac{1}{k+1} \left(\left( \frac{\mu(II_0)}{\mu(I)}- \frac{\mu(I'I_0)}{\mu(I')}\right)a_n^k + \left( \frac{\mu(II_1)}{\mu(I)}- \frac{\mu(I'I_1)}{\mu(I')}\right)b_n^k\right)\right|\le\cr
&\le & \frac{|h(p_n)-h(q_n)|}{k+1}+\left|\frac{1}{k+1} \left( \frac{\mu(II_0)}{\mu(I)}- \frac{\mu(I'I_0)}{\mu(I')}\right)(a_n^k - b_n^k)\right|\cr
&\le &\frac{ (1-|p_n-q_n|)\log 2}{k+1}+|p_n-q_n|\frac{|a_n^k - b_n^k|}{k}\frac{k}{k+1}
\end{eqnarray}
We can rewrite  relation (\ref{stage1}) in the following way
$$\frac{|a_{n-1}^{k+1}-b_{n-1}^{k+1}|}{k+1}\le \frac{(1-|p_n-q_n|)\log 2}{k+1}+|p_n-q_n|\frac{|a_n^k - b_n^k|}{k}\left(1-\frac{1}{k+1}\right)$$ and thus,
\begin{equation}\label{keyremark}
\Delta_{n-1}^{k+1}\le \frac{(1- |p_n-q_n|)\log 2}{k+1}+|p_n-q_n| \left(1-\frac{1}{k+1}\right)\Delta_n^k
\end{equation}
By lemma \ref{MA401} we then obtain,
\begin{equation}\label{keyremark2}
\Delta_{n-1}^{k+1}\le \max\left\{ \left(1-\frac{1}{k+1}\right)\Delta_n^k,\frac{ \log 2}{k+1}\right\}.
\end{equation} 
We use a recursion argument to finish the proof the lemma. First
observe that if for some $\ell\in\{0,...,k\}$ we have 
\begin{equation}\label{imaginary}
\Delta_{n+\ell}^{k-\ell}<\frac{\log 2}{k-\ell}
\end{equation} then we will also have $\ds \Delta_{n+\ell-1}^{k-\ell+1}<\frac{\log 2}{k-\ell+1}$, by relation (\ref{keyremark2}), and therefore $\ds\Delta_{n-1}^{k+1}\le\frac{\log 2}{k+1}$.

On the other hand, if inequality (\ref{imaginary}) does not hold for any $\ell\in\{0,...,k\}$ then by (\ref{keyremark2}) we get 
$$\Delta_{n+\ell-1}^{k-\ell+1}\le  \left(1-\frac{1}{k-\ell+1}\right)\Delta_{n+\ell}^{k-\ell}$$
and finally
\begin{equation}\label{finallemma}
\Delta_{n-1}^{k+1}\le \prod_{\ell=1}^{k+1}\left(1-\frac{1}{\ell+1}\right)\log 2\le \frac{e^2\log 2}{k+1}
\end{equation}
Take $\ds \eta(k)=\frac{e^2\log 2}{k+1}$. By equations (\ref{imaginary}) and (\ref{finallemma}) we get 
$\ds \Delta_{n-1}^{k+1}\le \max\left\{ \frac{\log 2}{k+1}, \frac{e^2\log 2}{k}\right\}=\frac{e^2\log 2}{k+1}=\eta(k)$ and the proof is complete.
\end{Proof}
We will also use the following two theorems of \cite{BH} that we include without proof for the convenience of the reader (a straight forward proof -without use of these theorems- is possible but much longer). 

\begin{Theorem}\label{t1} \cite{BH} Let $m$ be a probability measure in $[0,1)^D$ equipped with the filtration of $\ell$-adic cubes, $\ell\in\N$. Then
$$\dim_* (m)\leq h_*(m)\ .$$
Moreover, the following properties are equivalent :
\begin{enumerate}
\item $\dim_* (m)=h_*(m)$
\item$\dim_* (m)=\dim^*(m)=h_*(m)$
\item There exists a subsequence $(n_k)_{k\in\N}$ such that 
for $m$-almost every $x\in [0,1)^D$,
$$\lim_{k\rightarrow +\infty}\frac{\log m(I_{n_k}(x))}{-n_k\, \log \ell}
=\dim_* (m)\ .$$
\end{enumerate}
\end{Theorem}
\begin{Theorem}\label{t2} \cite{BH}  We also have
{\em $$h^*(m)\leq\Dim^* (m) ,$$}
and the following properties are equivalent :
\begin{enumerate}
\item {\em $\Dim^* (m)=h^*(m)$}
\item {\em $\Dim_* (m)=\Dim^* (m)=h^*(m)$}
\item There exists a subsequence $(n_k)_{k\in\N}$ such that 
for $m$-almost every $x\in [0,1)^D$,
{\em $$\lim_{k\rightarrow +\infty}\frac{\log m(I_{n_k}(x))}{-n_k\, \log \ell}
=\Dim^* (m)\ .$$}
\end{enumerate}
\end{Theorem}
\section{Proofs of the theorems}
To prove theorem \ref{letsgo1} we will use the following strong law of large numbers (cf. \cite{HH}).
\begin{Theorem}\label{Grandsnombres}({\bf Law of Large Numbers})
Let $(X_n)_{n\in\N}$ be a sequence of uniformly bounded in ${\mathcal L}^2$ real random variables on a
probability space 
$({\mathbb X},{\mathcal B}, P)$ and let $(\Im_n)_{n\in\N}$ be an increasing
sequence of $\sigma$-subalgebras of $\mathbb B$ such that $X_n$ 
is measurable with respect to $\Im_n$ for all
$n\in\N$. Then
\begin{equation}\label{grandsnombres} 
\ds \lim_{n\to\infty}\frac{1}{n}\sum_{k=1}^n
\left(X_k-\Esp (X_k|\Im_{k-1})\right)=0\hp  \hp P\mbox{-almost
surely}
\end{equation}
\end{Theorem}
Remark that the assumptions on the random variables are not optimal but it will be sufficient for our goal. The space here is $\K$, the filtration will be the dyadic one and $\mu$ will take the place of the probability measure $P$.
\begin{Proof}{\bf of Theorem \ref{letsgo1}.}
Consider the random variables $X_n$, $n\in\N$, defined on $\K$, given by
$$X_n(x)=\log\frac{\mu\left(I_n(x)\right)}{\mu\left(I_{n-1}(x)\right)},$$
where, for $x\in\K$, we have noted $I_n(x)$ the unique element of $\Im_n$ containing $x$.
The previous lemma implies that for all positive $p$'s
\begin{equation}\label{grandsnombres2}
\lim_{n\to\infty} \frac{1}{(n+1)}\sum_{j=1}^n
\left( \frac{1}{p}\sum_{k=1}^p\left[ X_{jp+k}- \Esp (X_{jp+k}|\Im_{jp})\right]\right)=0\hp,  \hp \mu\mbox{-almost
surely.}
\end{equation}
On the other hand, on each $I\in \Im_n$, the conditional expectation    
$\displaystyle \frac{1}{p}\sum_{k=1}^p\Esp (X_{np+k}|\Im_{np})$ 
is given by 
\begin{equation}\label{final0}
\frac{1}{p}\sum_{k=1}^p\Esp (X_{np+k}|\Im_{np})= \frac{1}{p}\sum_{K\in{\mathcal F}_p}\log\left(\frac{\mu(IK)}{\mu(I)}\right)\frac{\mu(IK)}{\mu(I)}.
\end{equation}
By proposition \ref{mainpro}, for every $\epsilon>0$ there exists $p\in\N$ such that for all $n\in\N$ and all $I$ in $\Im_n$
\begin{equation}\label{final1}
\left|\frac{1}{p}\sum_{K\in{\mathcal F}_p}\log\left(\frac{\mu(IK)}{\mu(I)}\right)\frac{\mu(IK)}{\mu(I)}-c_{n}\right|<\epsilon,
\end{equation}
where $c_{n}$ is a constant depending only on $n$ and on the chosen $p$ but not on the cylinder $I$ of $\Im_n$.

It is also easy to see that the variable $(X_n)_{n\in\N}$ are uniformly bounded in ${\mathcal L}^2(\mu)$. 
We deduce, using the relations (\ref{grandsnombres2}) and (\ref{final0}), that for every $\epsilon>0$ there exists $p\in\N$ and a sequence $(c_n)_{n\in\N}$ of real numbers such that
\begin{eqnarray}\label{ultime}
-\epsilon &<& \liminf_{n\to\infty} \frac{1}{(n+1)} \sum_{j=1}^n
\left( \frac{1}{p}\sum_{k=1}^pX_{jp+k}- c_j\right)\leq \cr
&\leq &\limsup_{n\to\infty} \frac{1}{(n+1)}\sum_{j=1}^n
\left( \frac{1}{p}\sum_{k=1}^pX_{jp+k}- c_j\right)<\epsilon,
\end{eqnarray}
$\mu$-almost everywhere on $\K$.
This relation implies that  
\begin{equation}\label{final2}
\liminf_{n\to\infty} \frac{-1}{(n+1)}\sum_{j=1}^n
 c_j -\epsilon < \liminf_{n\to\infty} \frac{-1}{p}\frac{1}{(n+1)}\sum_{j=1}^n \sum_{k=1}^pX_{jp+k} < 
\liminf_{n\to\infty} \frac{-1}{(n+1)} \sum_{j=1}^n c_j+\epsilon
\end{equation} and 
\begin{equation}\label{final2bis}
\limsup_{n\to\infty} \frac{-1}{(n+1)}\sum_{j=1}^n
 c_j -\epsilon < \limsup_{n\to\infty} \frac{-1}{p}\frac{1}{(n+1)}\sum_{j=1}^n \sum_{k=1}^pX_{jp+k} < 
\limsup_{n\to\infty} \frac{-1}{(n+1)} \sum_{j=1}^n c_j+\epsilon
\end{equation} 
$\mu$-almost everywhere on $\K$.
If we note 
$$ \underline{c}=\liminf_{n\to\infty} \frac{-1}{(n+1)\log 2}\sum_{j=1}^nc_j \mbox{ and }\overline{c}=\limsup_{n\to\infty} 
\frac{-1}{(n+1)\log2}\sum_{j=1}^n c_j,$$ 
we deduce from (\ref{final2}) and (\ref{final2bis})  that $\dim_*\mu=\underline{c}$ and $\Dim_*\mu=\overline{c}$.

Furthermore, the inequalities (\ref{ultime}) imply that for every positive $\epsilon$ there is a strictly increasing sequence of natural numbers $(n_l)_{l\in\N}$ verifying 
$$-\epsilon<\liminf_{l\to\infty} \frac{-1}{(n_l+1)} \sum_{j=1}^{n_l}
\left( \frac{1}{p}\sum_{k=1}^pX_{jp+k}\right)-\underline{c} \le \limsup_{l\to\infty} \frac{-1}{(n_l+1)}\sum_{j=1}^{n_l}
\left( \frac{1}{p}\sum_{k=1}^pX_{jp+k}\right)- \underline{c}<\epsilon.$$
One easily proves (using, for instance, Cantor's diagonal argument) that there exists a strictly increasing sequence of natural numbers $(n_l)_{l\in\N}$ such that
$$\lim_{l\to\infty}\frac{-1}{n_l\log 2}\log\mu\left(I_{n_l}(x)\right)=\dim_*(\mu),$$ for $\mu$-almost all $x\in\K$. 

Similarly, there exists a strictly increasing sequence of natural numbers $(\hat n_l)_{l\in\N}$ such that
$$\lim_{l\to\infty}\frac{-1}{\hat n_l\log 2}\log\mu\left(I_{\hat n_l}(x)\right)=\Dim_*(\mu),$$ for $\mu$-almost all $x\in\K$. We use theorems \ref{t1} and \ref{t2} to finish the proof.
\end{Proof}

To prove theorem \ref{letsgo2} we will use proposition \ref{mainpro} and lemma \ref{Grandsnombres}.
\begin{Proof}{\bf of theorem \ref{letsgo2}}
Take $\epsilon>0$ and let $(p_n,q_n)_{n\in\N}$ and $(p_n',q_n')_{n\in\N}$ be two sequences of weights satisying $0<p_n,q_n,p_n',q_n'<1$ for all $n\in\N$ and 
$$||(p_n,q_n)_{n\in\N}-(p_n',q_n')_{n\in\N}||_{\infty} <\zeta.$$ We note $\mu$ and $\mu'$ the measures corresponding to these two sequences of weights. We will show that 
$$|\dim_*(\mu)-\dim_*(\mu')|<\epsilon,$$
if $\zeta$ is small enough.

It follows from proposition \ref{mainpro} that there exist  a natural number $p$ large enough and two sequences of real numbers $(c_n)_{n\in\N},(c_n')_{n\in\N}$ such that the following relations hold:
$$\left|\frac{1}{p}\sum_{K\in{\mathcal F}_p}\log\left(\frac{\mu(IK)}{\mu(I)}\right)\frac{\mu(IK)}{\mu(I)}-c_{n}\right|<\frac{\epsilon}{4}$$ and 
$$\left|\frac{1}{p}\sum_{K\in{\mathcal F}_p}\log\left(\frac{\mu'(IK)}{\mu'(I)}\right)\frac{\mu'(IK)}{\mu'(I)}-c_{n}'\right|<\frac{\epsilon}{4}$$ 
for all cylinders $I\in \Im_{np}$ and all $n\in \N$.
Since $p$ is a fixed finite number it suffices to take $\zeta$ small in order to have 
$$\left|\frac{1}{p}\sum_{K\in{\mathcal F}_p}\log\left(\frac{\mu(IK)}{\mu(I)}\right)\frac{\mu(IK)}{\mu(I)}-\frac{1}{p}\sum_{K\in{\mathcal F}_p}\log\left(\frac{\mu'(IK)}{\mu'(I)}\right)\frac{\mu'(IK)}{\mu'(I)}\right|<\frac{\epsilon}{2},$$
for all $I\in \Im_{np}$ and all $n\in \N$.
Hence, 
$$-\epsilon <\liminf_{n\to\infty} \frac{1}{(n+1)} \sum_{j=1}^n |c_j-c_j'|\leq \limsup_{n\to\infty} \frac{1}{(n+1)}\sum_{j=1}^n|c_j-c_j'|<\epsilon.$$
we deduce from (\ref{final2}) and (\ref{final2bis})  that 
$ |\dim_*(\mu)-\dim_*(\mu')|<\epsilon$ and $|\Dim_*(\mu)-\Dim_*(\mu')|<\epsilon, $ which completes the proof.
\end{Proof}
Theorem \ref{letsgo2} have a limited validity as we show in the following section. 

\section{A counterexample}
For every $\epsilon>0$ we construct two dyadic doubling measures $\mu$ and $\nu$ on $\K$ such that if $\displaystyle X_n(x)=\log\frac{\mu\left(I_n(x)\right)}{\mu\left(I_{n-1}(x)\right)}$ and $\displaystyle Y_n(x)=\log\frac{\nu\left(I_n(x)\right)}{\nu\left(I_{n-1}(x)\right)},$ for $n\in\N$ then 
\begin{equation}\label{condchercheeCantor}
\sup_{n\in\N} \Big|\Big| X_n-Y_n \Big|\Big|_{L^\infty}<\epsilon 
\end{equation}
 and 
$\displaystyle\dim_*(\mu)-\dim_*(\nu)|>\frac{1}{4}$.
A first example was proposed to us by professor Alano Ancona; the proof provided here is of a similar nature.

The construction is carried out in two stages. We start by finding two Bernoulli measures  satisfying \ref{condchercheeCantor}
and afterwards we will modify them by a reccurent process to get the corresponding dimensions very different.

For $I\in\Im_n$ we note $\widehat I$ the unique cylinder of the 
$(n-1)$th generation $\Im_{n-1}$ contenant $I$. Relation (\ref{condchercheeCantor}) can now be reformulated in the following way 
\begin{equation}\label{condchercheeCantor2}
\ds\left|\frac{\mu(I)}{\mu(\widehat I)} : 
\frac{\nu(I)}{\nu(\widehat {I})}-1 \right|<\epsilon  \hp, \hp
\mbox{ for all cylinders }I \mbox{ of } \bigcup_{n\in\N}\Im_n
\end{equation}

\subsubsection{The starting point}
Take $\epsilon>0$ and  $\lambda_0$ the Lebesgue (uniform) measure (of dimension 1) on $\K$.

Consider the Bernoulli measure $\rho_0$  of weight variable $\frac{1}{2}-\epsilon$, i.e. such that for $I\in\Im_n$, $n\in\N$,

\begin{equation}\label{rho-0}
\rho_0(II_0)=(\frac{1}{2}-\epsilon)\rho(I)\hp , \hp
\rho_0(II_1)=(\frac{1}{2}+\epsilon)\rho(I). 
\end{equation}

Put $\mu_0=\lambda_0$ and $\nu_0=\rho_0$.
By construction the measures  $\lambda_0$ and $\rho_0$ verify  condition
(\ref{condchercheeCantor}), are exact and doubling on the dyadics.
Moreover, we have
$$\dim\rho_0=h_*(\rho)= -\frac{\left(\frac{1}{2}-\epsilon\right)}{\log 2}\log\left(\frac{1}{2}-\epsilon\right)-\frac{\left(\frac{1}{2}+\epsilon\right)}{\log 2}\log\left(\frac{1}{2}+\epsilon\right)$$ 
It is clear that $\lambda_0$ and $\rho_0$ are singular.
Furthermore by the Shannon-MacMillan formula (cf for instance \cite{Zin}),
$$\lim_{n\to\infty}\frac{\log\rho_0(I_n(x))}{n}=h_*(\rho_0)\hp
\hp\rho_0\mbox{-almost evcerywhere on} \K.$$

Hence, we can find $n_1\in\N$ and a partition $\{F_0,F_1\}$ of
$\Im_{n_1}$ verifying :

\begin{itemize}
\item[1.] $F_0\cup F_1=\Im_{n_1}$
\item[2.]  $\left| \ds \frac{\log\rho_0(I)}{n} +h_*(\rho_0)\right|
< \epsilon$ for all $I\in F_1$
\item[3.]  $\left| \ds \frac{\log\lambda_0(I)}{n} +\log 2\right|
< \epsilon$ for all $I\in F_0$
\item[4.] $\ds\sum_{I\in F_1}\rho_0(I)>1-\epsilon$
\item[5.] $\ds\sum_{I\in F_0}\lambda_0(I)>1-\epsilon$
\end{itemize}
 Let us also define the Bernoulli measures $\lambda_1$
and $\rho_1$ on $\K$ in the following way 
\begin{eqnarray}\label{lambda1-rho1}
\rho_1(I_0)=\delta\nonumber & \mbox{ and }&  \rho_1(I_1)=1-\delta \nonumber
\\*[-2.5mm]
\\*[-2.5mm]
\lambda_1(I_0)=\delta(1-\epsilon) &\mbox{ and}&
\hp \lambda_1(I_1)=1-\delta(1-\epsilon)\nonumber
\end{eqnarray}
where $\delta>0$ will be fixed later.

\subsubsection{Going on with the construction}
For $I_{i_1...i_n}\subset I\in F_1$ we put
\begin{eqnarray}
\mu_1(I_{i_1...i_n})&=&\mu_0(I_{i_1...i_{n_1}})\lambda_1(I_{i_{n_1}...i_n})
\nonumber \\*[-2.5mm]
\\*[-2.5mm]
\nu_1(I_{i_1...i_n})&=&\nu_0(I_{i_1...i_{n_1}})\rho_1(I_{i_{n_1}...i_n})
 \nonumber
\end{eqnarray}
and  for $I_{i_1...i_n}\subset I\in F_0$ 
\begin{equation}
\mu_1(I_{i_1...i_n})=\mu_0(I_{i_1...i_{n}})\hph
\mbox{ , }\hph
\nu_1(I_{i_1...i_n})=\nu_0(I_{i_1...i_{n}})
\end{equation}

Remark that for $I=I_{i_1...i_n}$ with $n\le n_1$ we leave
$\mu_1(I)=\mu_0(I) $ and $\nu_1(I)=\nu_0(I)$.

The restrictions of the measures $\mu_1$ and $\nu_1$ in the cylinders of
$\Im_{n_1}=F_0\cup F_1$ are Bernoulli measures of different dimensions, so they are singulars between them. 
Therefore, we can find $n_2\in\N$ and a partition $\{F_{00}, F_{01}, F_{10},F_{11}\}$ of
$\Im_{n_2}$  such that
\begin{itemize}
\item[1.] $I\in F_{j0}\cup F_{j1}$ if and only if there is $J\in F_j$ such that $I\subset J$, $j\in\{0,1\}$.
\item[2.] $\ds\left|\frac{\log\mu_1(I)}{n_2}+\log 2\right|<\epsilon^2$
for all $I\in F_{00}$.
\item[3.]  $\ds\left|\frac{\log\nu_1(I)}{n_2}+h_*(\rho_1)\right|<\epsilon^2$
for all $I\in F_{11}$.
\item[4.] $\ds\sum_{\stackrel{\scriptstyle J\in F_{00}}{J\subset
I}}\mu_1(J)>(1-\epsilon^2)\mu_1(I)$\hp  et\hp  $\ds\sum_{\stackrel{\scriptstyle
J\in F_{01}}{J\subset 
I}}\nu_1(J)>(1-\epsilon^2)\nu_1(I)$ pour  $I\in F_0$ 
\item[5.]  $\ds\sum_{\stackrel{\scriptstyle J\in F_{10}}{J\subset
I}}\mu_1(J)>(1-\epsilon^2)\mu_1(I)$\hph and \hph $\ds\sum_{\stackrel{\scriptstyle
J\in F_{11}}{J\subset 
I}}\nu_1(J)>(1-\epsilon^2)\nu_1(I)$ pour  $I\in F_1$ 
\end{itemize}

If $I\in F_{00}\cup F_{10}$ and  $\ds
J\in\bigcup_{n\in\N}\Im_n$, we put 
$$\mu_2(IJ)=\mu_1(I)\lambda_0(J)\hph\mbox{ ,
}\hph\nu_2(IJ)=\nu_1(I)\rho_0(J).$$ 
If   $I\in F_{01}\cup F_{11}$ and  $\ds J\in\bigcup_{n\in\N}\Im_n$ 
we put
$$\mu_2(IJ)=\mu_1(I)\lambda_1(J)\hph\mbox{ ,
}\hph\nu_2(IJ)=\nu_1(I)\rho_1(J).$$ 
Finally, for  $I\in\Im_n$, with $n\le n_2$ we keep the same mass distribution $\mu_2(I)=\mu_1(I)$ et
$\nu_2(I)=\nu_1(I)$. 

Suppose the measures $\mu_k$, $\nu_k$ and the partition $\bigl\{F_{i_1...i_k}\hp,\hp
i_1,...,i_k\in\{0,1\}\bigr\}$ de $\Im_{n_k}$ are constructed.
As in the two first stages, the restrictions of the measures 
$\mu_k$ and $\nu_k$ on such cylinder of 
$\Im_{n_k}$ are supposed to be Bernoulli measures: whether 
$\lambda_0$ and $\rho_0$ whether $\lambda_1$ and$\rho_1$, respectively.

The measures $\mu_k$ and $\nu_k$ are singular between them. Hence, there is  $n_{k+1}>n_k$ and a partition 
$\bigl\{F_{i_1...i_{k+1}}\hp,\hp  
i_1,...,i_{k+1}\in\{0,1\}\bigr\}$ of $\Im_{n_{k+1}}$ satisfying

\begin{itemize}
\item[1.]  $I\in F_{i_1...i_k0}\cup F_{i_1...i_k1}$ if and only if there is $J\in F_{i_1...i_k}$ such 
that $I\subset J$, with $i_1,...,i_k\in\{0,1\}$
\item[2.]
$\ds\left|\frac{\log\mu_k(I)}{n_{k+1}}+\log 2\right|<\epsilon^{k+1}$
for all $I\in F_{i_1...i_{k-1}00}$.
\item[3.]  $\ds\left|\frac{\log\nu_k(I)}{n_2}+h_*(\rho_1)\right|<\epsilon^{k+1}$
for all $I\in F_{i_1...i_{k-1}11}$.
\item[4.] $\ds\sum_{\stackrel{\scriptstyle J\in F_{i_1...i_{k-1}00}}{J\subset
I}}\mu_k(J)>(1-\epsilon^{k+1})\mu_k(I)$\hph  and\hph
$\ds\sum_{\stackrel{\scriptstyle 
J\in F_{i_1...i_{k-1}01}}{J\subset 
I}}\nu_k(J)>(1-\epsilon^{k+1})\nu_k(I)$,\\ 
for all  cylinders $I\in F_{i_1...i_{k-1}0}$. 
\item[5.]  $\ds\sum_{\stackrel{\scriptstyle J\in F_{i_1...i_{k-1}10}}{J\subset
I}}\mu_k(J)>(1-\epsilon^{k+1})\mu_k(I)$\hph et\hph
 $\ds\sum_{\stackrel{\scriptstyle
J\in F_{i_1...i_{k-1}11}}{J\subset 
I}}\nu_k(J)>(1-\epsilon^{k+1})\nu_k(I)$,\\
 for all cylinders  $I\in F_{i_1...i_{k-1}1}$. 
\end{itemize}

If $I\in F_{i_1...i_k0}$,  $i_1,...,i_k\in\{0,1\}$, then for all
$\ds J\in\bigcup_{n\in\N}\Im_n $  we put
$$\mu_{k+1}(IJ)=\mu_k(I)\lambda_0(J)\hph\mbox{ and}\hph
\nu_{k+1}(IJ)=\nu_k(I)\rho_0(J).$$ 
If $I\in F_{i_1...i_k1}$, $i_1,...,i_k\in\{0,1\}$, then for all
$\ds J\in\bigcup_{n\in\N}\Im_n $  we put
$$\mu_{k+1}(IJ)=\mu_k(I)\lambda_1(J)\hph\mbox{ and}\hph
\nu_{k+1}(IJ)=\nu_k(I)\rho_1(J).$$ 

\subsubsection {Properties of the measures defined}

It is clear that the sequences  $(\mu_n)_{n\in\N}$ and
$(\nu_n)_{n\in\N}$  converge towards two probability measures $\mu$
and $\nu$ respectively.
By the construction $\mu$ and $\nu$ are doubling on the dyadics,
exacts and  satisfy 
(\ref{condchercheeCantor}). 

On the other hand, clearly $\dim_*\mu=1$ and it is not difficult to see that
$\dim_*\nu\le\frac{1}{2}$, if $\delta $ is small enough, since 
$\ds\liminf_{n\to\infty}\frac{-\log \nu(I_n(x))}{n\log 2}=\frac{h_*(\rho_1)}{\log 2}$, $\nu$-almost everywhere.
Evenmore, the measures $\mu$ and $\nu$ satisfy the conclusion of theorem \ref{letsgo1}. The counterexample is complete.

\begin{Aknow} The author would like to thank pr. M. Babillot for having carefully read this work and for her numerous comments.
\end{Aknow}
\bibliographystyle{alpha}
\bibliography{biblio.bib}
\vspace{1cm}
\hfill
\begin{tabular}{c}
{\bf Athanasios BATAKIS}\\
{MAPMO}\\
{Universit\'e d'Orl\'eans}\\
{BP 6759}\\
{45067 Orl\'eans cedex 2}\\
{FRANCE}\\
{email: batakis@labomath.univ-orleans.fr}
\end{tabular}
\end{document}